\documentclass[11pt]{article}
\usepackage{amsmath, amsfonts, amssymb, amsthm}
\usepackage[square]{natbib}

\usepackage{newejpart}

\usepackage{float}
\usepackage{amsmath}
\usepackage{graphicx}
\usepackage{latexsym}
\usepackage{amsfonts}
\usepackage{amssymb}%

\usepackage{xspace}

\newcommand{\eps}{\varepsilon}

\numberwithin{equation}{section}

\renewcommand{\th}{\ensuremath{^{\text{th}}}\xspace}

\newtheorem{theorem}{Theorem}[section]
\newtheorem{corollary}[theorem]{Corollary}
\newtheorem{lemma}[theorem]{Lemma}

\theoremstyle{definition}
\newtheorem{remark}[theorem]{Remark}

\newcommand{\E}{{\mathbb E}}
\newcommand{\R}{{\mathbb R}}

\newcommand{\kappah}{{\hat{\kappa}}}
\newcommand{\nt}{\tilde{n}}
\newcommand{\vm}{\vec{m}}
\newcommand{\vM}{\vec{M}}
\newcommand{\vx}{\vec{x}}
\newcommand{\vy}{\vec{y}}
\newcommand{\vp}{\vec{p}}
\newcommand{\vxi}{\vec{\xi}}
\newcommand{\vg}{\vec{\eta}}
\newcommand{\vc}{\vec{x}}
\newcommand{\lhe}{\stackrel{{\displaystyle <}}{\mbox{\scriptsize\_\,\_}}}
\newcommand{\wto}{\stackrel{\mbox{\scriptsize \!\!\rm w}}{\to}}
\newcommand{\apxopt}{\tilde{\alpha}}

\title{Two-player Knock~'em Down}
\amsprimary{91A60}
\amssecondary{91A05}
\startpage{198}
\endpage{212}
\submissiondate{December 14, 2006}
\acceptancedate{January 31, 2008}
\papernumber{9}
\ejpvolume{13}
\ejpyear{2008}
\abstract{We analyze the two-player game of Knock~'em Down, asymptotically as
  the number of tokens to be knocked down becomes large.  Optimal play
  requires mixed strategies with deviations of order $\sqrt{n}$ from
  the na\"ive law-of-large numbers allocation.  Upon rescaling by
  $\sqrt{n}$ and sending $n\to\infty$, we show that optimal play's
  random deviations always have bounded support and have marginal
  distributions that are absolutely continuous with respect to
  Lebesgue measure}
\author{\large James Allen Fill\thanks{The research of J.A.F.\ was supported by NSF grants DMS-0104167 and DMS-0406104 and by The Johns Hopkins University's Acheson~J.~Duncan Fund for the Advancement of Research in Statistics.  Research carried out in part while a visiting researcher at Microsoft Research.}\\[-2pt]
\small Department of Applied Mathematics and Statistics\\[-2pt]
\small The Johns Hopkins University\\[-2pt]
\small Baltimore, Maryland 21218-2682\\[-2pt]
\small U.S.A. \\
\small \url{http://www.ams.jhu.edu/~fill/} 
\and \large David B. Wilson \\[-2pt]
 \small Microsoft Research\\[-2pt]
 \small One Microsoft Way\\[-2pt]
 \small Redmond, Washington 98052\\[-2pt]
 \small U.S.A.\\
 \small \url{http://dbwilson.com}
}

\begin{document}

\maketitle

\section{Introduction}
\label{S:intro}

\subsection{Knock~'em Down}
\label{S:KED}

In the game of Knock~'em Down, a player is given $n$~tokens which
(s)he arranges into~$k$ piles, or bins.  After that, a $k$-sided die
is thrown; if the outcome is side~$i$ (which occurs with probability
$p_i > 0$), then a token is knocked down from the~$i$\th pile.  In the
event that there are no tokens in the~$i$\th bin, then no tokens get
knocked down.  The die is thrown repeatedly until all the tokens have
been knocked down.  [In the original version of this game as described
in~\cite{Hunt}~\cite{BF_UMAP}, $n = 12$ and two fair six-sided dice
are thrown, with the bin chosen being given by ($1$~less, let us say,
than) the sum of the two numbers showing. In that case, $k = 11$ and
$p_i \equiv (6 - |i - 6|) / 36$.  The reader may wish to note that we
have altered the spelling from ``Knock~'m Down'' to ``Knock~'em Down''.]

We consider two versions of Knock~'em Down.  In solitaire Knock~'em
Down, which we analyze in a separate article, there is one player, and
his goal is to minimize the expected number of iterations until all
the tokens have been knocked down.  Solitaire Knock~'em Down is also
equivalent to a two-player zero-sum game, where the payoff to the
winner is the expected number of extra die throws that the other
player requires to knock down all his tokens, and the goal is to
maximize the expected payoff.  In competitive Knock~'em Down, which we
analyze in this article, it is enough merely to win, and the amount by
which a player wins is irrelevant.  There are two players, who each
arrange their tokens into bins without seeing what the other player is
doing, and then the same die is used to knock down tokens from each of
the players' bins.  The winner is the player whose tokens get all
knocked down first (the outcome may be a tie, which we resolve with a
coin flip).  With competitive Knock~'em Down, as we shall see, there
is an interesting Nash equilibrium.

\begin{figure}[htbp]
\includegraphics[width=\textwidth]{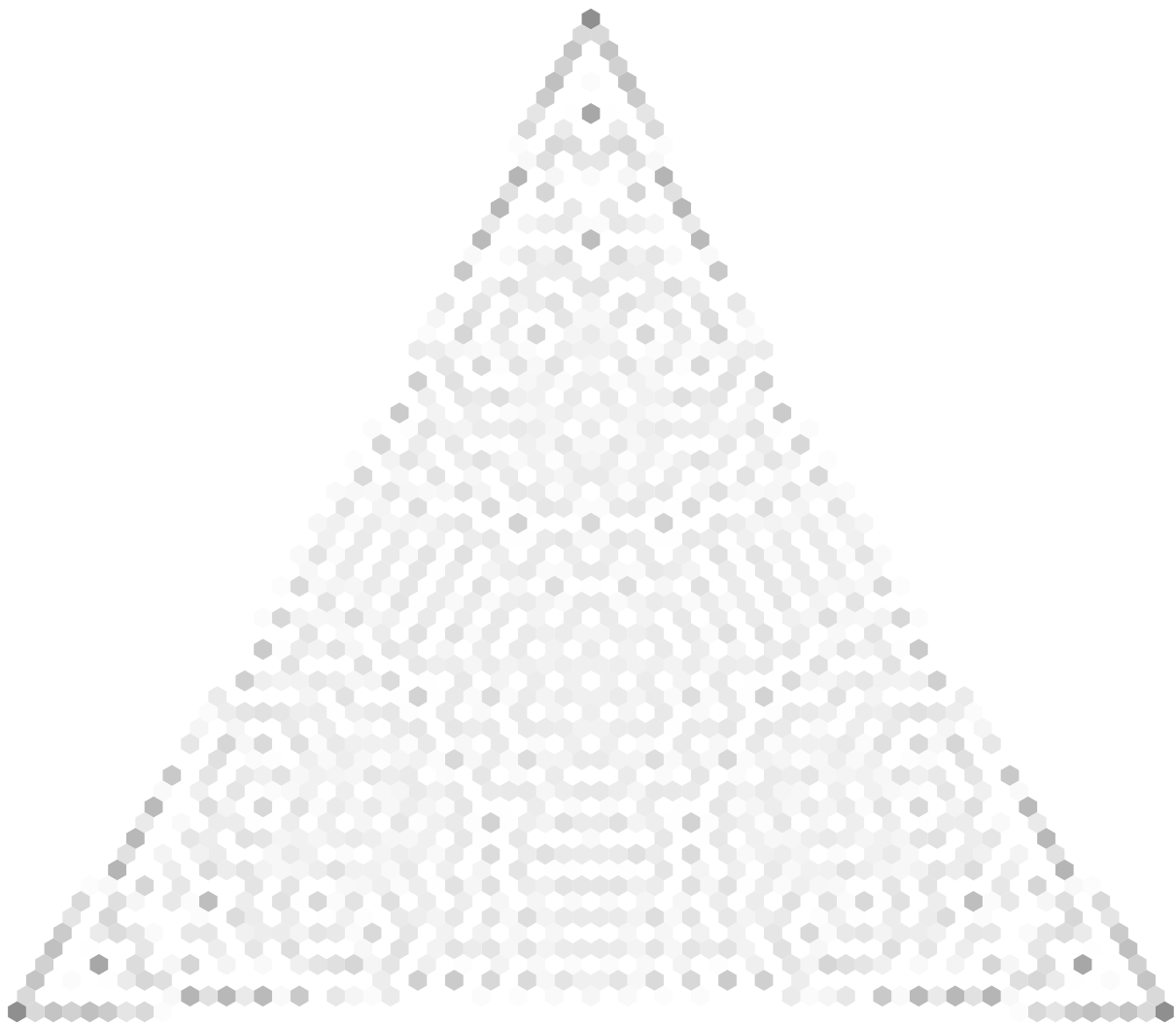}
\caption{
  An optimal strategy for (an approximation to) the continuous
  limiting version of two-player Knock~'em Down described in
  Section~\ref{S:sscont} when $k=3$ and $\vp=(1/3,1/3,1/3)$.  The
  strategy chooses a hexagon with probability proportional to its
  darkness, and then allocates the chips according to the hexagon's
  coordinates.  Hexagons near the lower-left, lower-right, or top of
  the triangle allocate more chips to the first, second, or third
  piles respectively.  The top corner has coordinates
  $(-0.28,-0.28,+0.56)/\sqrt3$.  }
\label{fig:150}
\end{figure}

Competitive Knock~'em Down is quite easy to analyze when $k = 2$.  The
result is Theorem~1 in~\cite{MR1767066}; the authors show that the
best strategy is to use the allocation $(m, n - m)$, where~$m$ is a
median of the Binomial$(n, p_1)$ distribution.  It is instructive to
consider competitive Knock~'em Down in the next simplest case, in
which a fair three-sided die is used.  The first player may guess that
the best strategy is to place $n/3$ tokens into each of the three bins
(assume for convenience that $n$ is divisible by $3$).  But if the
first player uses that strategy, then the second player can
``undercut'' by placing $(n/3) - 1$ tokens into each of the first two
bins and $(n/3) + 2$ tokens into the third bin.  Then, for large~$n$,
the probability that the second player's third bin empties out last is
only slightly larger than $1/3$, while the probability that the first
player's first or second bin finishes last is only slightly smaller
than $2/3$, so that the second player wins about two-thirds of the
time.  It turns out that an optimal strategy is to allocate
approximately $n/3$ tokens to each bin, but with certain random
perturbations to the bin allocations.  (See Figure~\ref{fig:150}.)  In
game-theoretic terminology, optimal play employs a mixed (non-pure)
strategy.

To simplify the analyses of both games for general~$k$ and~$\vp :=
(p_1, \dots, p_k)$, we suppose that the games are run in continuous
time, with the die being thrown at instants governed by a Poisson
point process with rate~$1$.  Throwing the die at random times rather
than at deterministic times has absolutely no impact on the outcome of
a game of competitive Knock~'em Down; likewise, in the case of
solitaire Knock~'em Down the expected (total) clearance time (i.e.,
time to knock down all tokens) remains unchanged.  Suppose in either
game that a player places~$\xi_i$ tokens into bin~$i$, and let~$T_i$
denote the time that it takes for bin~$i$ to be cleared. Since the
game is run in continuous time, the $T_i$'s are mutually independent,
and the variable $p_i T_i$ is the sum of~$\xi_i$ independent
exponential random variables with unit mean, i.e., $p_i T_i$ is a
Gamma random variable with shape parameter $\xi_i$.  The clearance
time is $T = \max_i T_i$.

\subsection{Asymptotic scale:\ root-$n$ deviations}
\label{S:root}

Partial results concerning optimal play have been derived by Art
Benjamin, Matt Fluet, and Mark
Huber~\cite{BF_UMAP}~\cite{MR1767066}~\cite{MR2002g:91048}~\cite{Fluet}
for both versions of Knock~'em Down, but in neither case has optimal
play been characterized analytically.  For both versions, we focus
here on the asymptotics of optimal play when the die remains fixed
[i.e.,\ the vector $\vp$ is held constant] and the number of
tokens~$n$ becomes large.  While optimal play in competitive Knock~'em
Down requires a mixed strategy, it is useful to consider first what
happens when a player deterministically places~$\xi_i$ tokens into the
$i$\th bin ($1\leq i\leq k$).  If~$\xi_i$ is large, then by the
central limit theorem, $T_i$ is approximately normal with mean $\xi_i
/ p_i$ and variance $\xi_i / p_i^2$.  In particular, bin $i$ is
emptied by time $T_i = (\xi_i / p_i) + O_p(\sqrt{\xi_i} / p_i)$.  (We
have used here the usual notation $X_n = O_p(y_n)$ to mean that
$\inf_n \Pr[|X_n| \leq c y_n]$ tends to~$1$ as $c \rightarrow \infty$;
equivalently, one says that the family of distributions
$\mathcal{L}(X_n / y_n)$ is \textit{tight}.)

If $\xi_i \equiv n q_i$, where the $q_i$'s are nonnegative and sum to
unity, then the clearance time~$T$ satisfies $T = n \max_i (q_i / p_i)
+ O_p(\sqrt{n})$.  Since $\max_i (q_i / p_i) \geq 1$, this suggests
that the optimal allocation is $\xi_i \approx n p_i$, and indeed for
the solitaire game this was first shown in~\cite{MR2002g:91048}.  We
show in our companion paper on solitaire Knock~'em Down that the
optimal choice is $\xi_i = n p_i + O(\sqrt{n})$ for a judicious choice
of the $O(\sqrt{n})$ perturbations.  For competitive Knock~'em Down,
our Theorem~\ref{T:first} below shows similarly that an optimal mixed
strategy will never choose~$\xi_i$ differing from~$n p_i$ by more than
order $\sqrt{n}$.

Define the overplay of an allocation~$\vxi$ relative to $n \vp$ to be
$$
\max_i \left( \frac{\xi_i}{p_i} - n \right).
$$
In the following lemma, $\{\mbox{$\vxi$ beats~$\vg$}\}$ is the event
that the corresponding clearance times satisfy the strict inequality
$T_{\vxi} < T_{\vg}$, or else $T_{\vxi} = T_{\vg}$ and $\vxi$ wins the
coin toss.  The notation extends naturally to
$\{\mbox{$\alpha$ beats~$\beta$}\}$ for mixed strategies $\alpha,
\beta$; see the start of Section~\ref{S:continuous} for careful
analogous definitions in the continuous-game setting that we introduce
in Section~\ref{S:sscont}.

\begin{lemma} \label{L:overplay}
  If the overplay of allocation~$\vg$ exceeds the overplay of
  allocation~$\vxi$ by at least $\left. w \sqrt{k n} \right/ \min_i p_i$,
  then $\Pr[\text{$\vg$ beats~$\vxi$}] < 1/w^2$.
\end{lemma}

\begin{proof}
  Let $j$ be the bin maximizing $\eta_j / p_j$.  In order for
  allocation $\vg$ to win, there must be some bin $\ell\neq j$ for
  which $\vg$'s bin $j$ clears out before $\vxi$'s bin $\ell$.  The
  clearance time $T_0$ of $\vg$'s bin $j$ has mean $\eta_j/p_j$ and
  variance $\eta_j/p_j^2\leq n/p_j^2$, while the clearance time $T_\ell$ of
  $\vxi$'s bin $\ell\neq j$ has mean $\xi_\ell/p_\ell \leq
  (\eta_j/p_j) - \left( \left. w \sqrt{k n} \right/ \min_i p_i
  \right)$ and variance $\leq \xi_\ell/p_\ell^2$.  Since $T_0$ and $T_\ell$
  are independent, by Chebyshev's inequality
  $$ \Pr[0<T_\ell-T_0] \leq \frac{(n/p_j^2) + (\xi_\ell/p_\ell^2)}{\big(w \sqrt{k n}/\min_i p_i\big)^2}\leq \frac{1+(\xi_\ell/n)}{w^2 k},$$
  and since the support of $T_\ell-T_0$ is $\R$, the first above inequality
  is strict.  Upon summing over $\ell\neq j$, we find that the
  probability that there is some bin $\ell\neq j$ for which $T_0 <
  T_\ell$ is less than $1/w^2$.
\end{proof}

\begin{theorem} \label{T:first}
  No optimal competitive Knock~'em Down (mixed) strategy ever
  overplays by more than $\left.\left(3\sqrt{3+\sqrt5} \sqrt{k n}
      + 1 \right) \right/ \min_i p_i$.
\end{theorem}
\begin{proof}
Let~$\alpha$ be any optimal strategy, and let $a(x)$ be the probability
that strategy~$\alpha$ picks an allocation with an overplay of~$x$ or more.
\enlargethispage{12pt}
We show successively, in each case using
Lemma~\ref{L:overplay}, that
$$
a\Bigg(\frac{w\sqrt{k n}+1}{\min_i p_i}\Bigg) < p_1,\quad\quad
  a\Bigg(\frac{2w\sqrt{k n}+1}{\min_i p_i}\Bigg) < p_2,\quad\quad
  a\Bigg(\frac{3w\sqrt{k n}+1}{\min_i p_i}\Bigg) = 0,
$$
for $w=\sqrt{3+\sqrt5}$, $p_1=\frac{-1+\sqrt5}{2}$, and $p_2=\frac{3-\sqrt5}{2}$.  Let $p=1-(1/w^2) = (1+\sqrt5)/4$.

Suppose first that $a((w\sqrt{k n}+1)/\min_i p_i) \geq p_1$.
Consider pure (i.e.,\ deterministic) strategy~$\beta$ which always
plays $n \vp$ (rounded to an integer vector summing to~$n$).  Because
of the rounding, $\beta$ may overplay $n \vp$ slightly, but never more
than $1/\min_i p_i$.  By Lemma~\ref{L:overplay}, strategy~$\beta$
beats~$\alpha$ with probability
$> p_1 \times p = 1/2$, contradicting the optimality of $\alpha$.

Suppose next that $a((2w\sqrt{k n}+1)/\min_i p_i) \geq p_2$.
With positive probability ($> 1-p_1$), strategy~$\alpha$ overplays by no more
than $(w\sqrt{k n} + 1) / \min_i p_i$, so we can define a
strategy~$\beta$ which plays according to strategy~$\alpha$
conditioned to overplay by no more than $(w\sqrt{k n}+1)/\min_i
p_i$.  When~$\beta$ plays against~$\alpha$, by Lemma~\ref{L:overplay},
$\beta$ wins with
probability $> (1-p_1) \times \frac{1}{2} + p_2 \times p=1/2$, which
is again a contradiction to the optimality of~$\alpha$.

Suppose finally that $a((3w\sqrt{k n}+1)/\min_i p_i) = \delta > 0$.
Let~$\beta$ be the strategy which plays according to~$\alpha$
conditioned to overplay by no more than $(3w\sqrt{k n}+1)/\min_i
p_i$.  When~$\beta$ plays against~$\alpha$, by Lemma~\ref{L:overplay},
strategy~$\beta$ wins with
probability $>[(1 - \delta) \times \frac12] + [\delta \times (1-p_2)
\times p] = \frac12$, which is again a contradiction.
\end{proof}

\newcommand{\Normal}{\operatorname{Normal}}
\newcommand{\Poisson}{\operatorname{Poisson}}
\newcommand{\GGamma}{\operatorname{Gamma}}

\subsection{A continuous game}
\label{S:sscont}

From Theorem~\ref{T:first} we see that any optimal strategy for
competitive Knock~'em Down plays allocations deviating by
$O_p(\sqrt{n})$ from the na\"ive law-of-large-numbers allocation $n
\vp$.  For given~$n$ and allocation~$\vxi$, it is thus natural to
define numbers~$x_i$ so that $\xi_i = p_i n + x_i \sqrt{n}$.  Then
\begin{align*}
T_i \sim& \GGamma(\xi_i,1)/p_i \\
    \dot{\sim}& \Normal\left(n+\frac{x_i}{p_i} \sqrt{n}, \frac{n}{p_i}
       + \frac{x_i}{p_i^2} \sqrt{n}\right) \\
     \dot{\sim}&
\,n + \Normal\left(\frac{x_i}{p_i} , \frac{1}{p_i}\right) \sqrt{n}.
\end{align*}
(By $\dot\sim$ we mean that the total variation distance between the
two distributions tends to $0$ as $n\to\infty$ with $p_i$ fixed and
$x_i$ bounded.)  The player chooses the $x_i$'s so that
$$ \sum_{i=1}^k x_i = 0$$
and $p_i n + x_i \sqrt{n}$ is an integer, and the player's tokens are
all exhausted at time approximately
$$ n + \sqrt{n} \max_i \left(\frac{x_i}{p_i}+\frac{Z_i}{\sqrt{p_i}}\right),$$
where the $Z_i$'s are independent standard normal random variables.

Thus the large-$n$ asymptotics of either version of Knock~'em Down is
effectively a continuous game (called solitaire/competitive continuous
Knock~'em Down), where the first player chooses real numbers~$x_i$
satisfying $\sum_i x_i = 0$, and his clearance time is
$$
T_{\vx} = \max_i \left(\frac{x_i}{p_i}+\frac{Z_i}{\sqrt{p_i}}\right).
$$
The second player similarly chooses numbers~$y_i$ summing to~$0$, with
clearance time~$T_{\vy}$ defined using the same~$Z_i$'s. The sequel
provides various rigorous connections between $n$-token Knock~'em Down
and our continuous game.

We present here two indications that even qualitative analysis of the
continuous game is not entirely trivial.  First, the na\"ive pure
strategy $\vc = \vec{0}$ has no optimal response.  Indeed, the
responding player can undercut by playing $(- \eps, - \eps, \dots, +
(k - 1) \eps)$, and letting $\eps \downarrow 0$ provides strategies
which give the respondent asymptotically the optimal probability $(k -
1) / k$ of winning, but no strategy achieves this probability.
Second, it is not immediately clear that our two-player continuous
game has a value (in the game-theoretic sense).  There are standard
tools for proving that a continuous game has a value, such as results
of Ky Fan~\cite{MR14:1109f}, but our payoff function is rather
severely discontinuous at certain points and the tools require a
semicontinuous payoff function.  Continuous games without values do
exist~\cite{MR20:265}, but it turns out that our game does indeed have
a value.  One can prove this by a suitable comparison of our game with
a ``ties go to player~1'' modification of the game having upper
semicontinuous payoff function, but we will give a somewhat more
direct proof whose basic idea is simply to pass to the limit from
$n$-token optimal strategies.

\subsection{Guide to later sections}
\label{guide}

With solitaire continuous Knock~'em Down, the optimal strategy is
deterministic, and we are able to characterize it for general $p_i$'s
\cite{FW:solitaire}.  With competitive continuous Knock~'em Down,
optimal play is random (i.e., mixed) with a rather complicated
distribution (see Figure~\ref{fig:150}).  Even in the simplest
nontrivial instance, where $k = 3$ and $\vp = (1 / 3, 1 / 3, 1 / 3)$,
we are unable to calculate an optimal strategy.  However, for
general~$k \geq 3$ and~$\vp$, we are able to derive some basic results
about optimal play.  In Section~\ref{S:continuous} we show, for
example, that any optimal strategy for the continuous game has
absolutely continuous marginals.  In Section~\ref{S:optexists} we
prove that the continuous game has an optimal strategy.  In
Section~\ref{S:asyopt} we show that a good strategy for the continuous
game can be converted to good strategies for the $n$-token games by
rounding; in particular, an optimal continuous strategy can be
converted to asymptotically optimal $n$-token strategies.  Finally, in
Section~\ref{S:open} we list some open problems arising from our work.

\section{Properties of optimal play}
\label{S:continuous}

The main result of this section (Theorem~\ref{T:ac}) asserts that the
marginals of any optimal strategy for the continuous game described in
Section~\ref{S:sscont} are absolutely continuous.  This will follow
from an estimate which holds (see Lemma~\ref{L:del}) equally well for
the discrete game.

We begin by establishing some terminology for the continuous game; analogous
terminology will be employed for the discrete game.  Recall that
$$
T_{\vx} = \max_i \left(\frac{x_i}{p_i}+\frac{Z_i}{\sqrt{p_i}}\right)
$$
is the clearance time for a player using allocation $\vx = (x_1,
\dots, x_k)$ with $x_1 + \cdots + x_k = 0$.  Here $Z_1, \dots, Z_k$
are independent standard normal random variables.  Recall that if
player~$1$ (say) uses allocation~$\vx$ and player~$2$ uses
allocation~$\vy$, then player~$1$ wins if and only if $T_{\vx} <
T_{\vy}$; for short, we simply say ``$\vx$ beats~$\vy$\,'', in which
case player~$2$ pays~$1$ unit (utile) to player~$1$.

The usual game-theoretic analysis takes the payoff $K(\vx, \vy)$ to be the
expected amount won by player~$1$, namely $\Pr[\mbox{$\vx$ beats~$\vy$}] -
\Pr[\mbox{$\vy$ beats~$\vx$}]$; for convenience we instead define 
$K(\vx, \vy)$ to be half this difference.

For any two mixed strategies~$\alpha$ and~$\beta$, $\Pr[\text{$\alpha$
  beats~$\beta$}]=\iint\!\Pr[\text{$\vx$ beats~$\vy$}]\,\alpha(d
\vx)\,\beta(d \vy)$ and we extend the definition of~$K$ to
\begin{align}
K(\alpha, \beta)
  &:= \iint\!K(\vx, \vy)\,\alpha(d \vx)\,\beta(d \vy)
    = \frac12\big(\Pr[\text{$\alpha$ beats~$\beta$}] - \Pr[\text{$\beta$ beats~$\alpha$}]\big)
\nonumber\\
\label{Kdef}
  &=  \frac{1}{2} \iint\!\left( \Pr[\text{$\vx$ beats~$\vy$}] -
        \Pr[\text{$\vy$ beats~$\vx$}] \right)\,\alpha(d \vx)\,\beta(d \vy).
\end{align}
We say that~$\alpha$ \emph{dominates}~$\beta$ if $K(\alpha, \beta) >
0$.  A strategy~$\alpha$ is \emph{optimal} if no mixed strategy
(equivalently, no pure strategy) dominates it.

Notice that payoffs are unaffected if we flip a fair coin to decide
ties.  That is, in the above definitions we may (and do) without
effect redefine $\Pr[\text{$\vx$ beats~$\vy$}]$ as $\Pr[T_{\vx} \lhe
T_{\vy}]$.  Here we have introduced the convenient notation $\Pr[A
\lhe B]$ to mean $\Pr[A < B] + \frac12 \Pr[A = B]$; similarly, $\Pr[A
\lhe B \lhe C]$ is shorthand for $\Pr[A<B<C] + \frac12 \Pr[A=B<C] +
\frac12 \Pr[A<B=C] + \frac16 \Pr[A=B=C]$.  (In our use of ``$\Pr[A
\lhe B \lhe C]$'' below, the third term always vanishes because
$\Pr[A<C]=1$.)

It will be convenient to view the contest between allocations~$\vx$
and~$\vy$ as follows.  Let $m_i := \max(x_i, y_i)$ for $1 \leq i \leq
k$.  For convenience we will refer to~$\vm = (m_1, \dots, m_k)$ as an
allocation, even though we are now working on the $\sqrt{n}$-scale of
deviations and the sum $m_1 + \cdots + m_k$ may exceed~$0$.
Correspondingly, on the scale of tokens, define $M_i := p_i n + m_i
\sqrt{n}$.  Let $I \equiv I_{\vm}$ be the bin that is cleared last
when the bin sizes are $M_1, \dots, M_k$.  [In the continuous game,
$I$ is defined correspondingly as $\arg\max_i (\frac{m_i}{p_i} +
\frac{Z_i}{\sqrt{p_i}})$, where the $Z_i$'s are independent standard
normals.]  If $x_I > y_I$, then~$\vy$ wins; if $x_I < y_I$, then $\vx$
wins; and if $x_I = y_I$, then the game is resolved by a coin flip.
(Thus, overall, $\Pr[\text{$\vx$ wins}]=\Pr[x_I\lhe y_I]$.)  To
compare various strategies, we will couple the $I_{\vm}$'s for various
values of~$\vm$, taking the viewpoint that the same random sequence of
die tosses will be made for a given game regardless of the
allocations~$\vx$ and~$\vy$ that the two players use (in the
continuous game, that the same vector of $Z_i$'s will be used).  Then
increasing~$m_i$ while leaving the other $m_j$'s ($j \neq i$) fixed
may change~$I$ from a value different from~$i$ to~$i$, but
otherwise~$I$ will not change.

\begin{lemma} \label{L:ichanges}
  If~$\vm$ is incremented by $1 / \sqrt{n}$ in position~$i$, then
  $$
  \Pr[\text{\rm $I$ changes}] \leq (1 + o(1)) / \sqrt{2 \pi p_i n} = 
  O(1 / \sqrt{n}).
  $$
\end{lemma}

\begin{proof}
  Let $M := M_i$, so that~$M$ is incremented by~$1$ to produce
  from~$\vm$ a new allocation~$\vm'$.  Ignore bin~$i$ for the moment,
  and let $T :=\max_{j \neq i} T_j$ be the time at which all remaining
  bins are cleared for allocation~$\vm$ (or~$\vm'$).  In order for
  $I_{\vm} \neq I_{\vm'}$, it must be that bin~$i$ was selected
  exactly~$M$ times up through time~$T$. Conditional on $T$, the
  number of times bin~$i$ was selected is Poisson distributed with
  parameter $\lambda := T p_i$, and the probability of exactly~$M$
  selections is
  $$
  \frac{e^{- \lambda} \lambda^M}{M!} \leq \frac{e^{- M} M^M}{M!} \leq 
  \frac{1}{\sqrt{2 \pi M}},
  $$
  Thus (unconditionally) $\Pr[\text{$I$ changes}]\leq 1/\sqrt{2\pi
    M}$.  We use this bound when $n$ is large and (say) $M \geq p_i n
  - n^{2 / 3}$.  When~$M$ is smaller, we instead use $\Pr[\text{$I$
    changes}] \leq \Pr[I_{\vm'} = i]$, which is $o(1 / \sqrt{n})$
  (indeed, is exponentially small in~$n$).
\end{proof}

\begin{lemma} \label{L:flat}
  For fixed~$\vp$, the distribution of $I_{\vm}$ is ``locally flat''
  as a function of~$\vm$: for bins~$h$, $i$, and~$j$ (not necessarily
  distinct), with~$e_h$ and~$e_i$ denoting unit vectors,
$$
\Pr[I_{\vm} = j]+\Pr[I_{\vm + \frac{e_h}{\sqrt{n}} + \frac{e_i}{\sqrt{n}}} = j]
 =\Pr[I_{\vm + \frac{e_h}{\sqrt{n}}} = j] + \Pr[I_{\vm + \frac{e_i}{\sqrt{n}}} = j] + O(1 / n).
 $$
\end{lemma}

\begin{proof}
  Let $F_t(\vM)$ denote the probability that the game with allocations
  given by~$\vM$ has clearance time $\leq t$.  If $Y_{\ell,t}$ denotes
  the number of times that bin~$\ell$ has been selected through
  time~$t$, then the $Y_{\ell}$'s are independent Poisson processes
  with respective intensity parameters~$p_{\ell}$.  We have
  $$
  F_t(\vM) = \Pr[\wedge_{\ell = 1}^k Y_{\ell, t} \geq M_\ell] =
   \prod_{\ell = 1}^k \Pr[\Poisson(p_{\ell} t) \geq M_{\ell}].
  $$
  The probability that the clearance time falls in $(t, t + dt)$ with
  bin~$j$ the last bin cleared is
$$
[F_t(\vM - e_j) - F_t(\vM)]\,p_j\,dt = \Delta_j F_t(\vM - e_j)\,p_j\,dt,
$$
where~$e_j$ is the~$j$\th unit vector and~$\Delta_j$ denotes the difference
operator
$$
\Delta_j f(\vx) = f(\vx) - f(\vx + e_j).
$$
Letting $\vM' = \vM - e_j$ we therefore have
$$
\Pr[I_{\vM} = j] = p_j \int_0^\infty\!\Delta_j F_t(\vM')\,dt,
$$
where here we are viewing $I = I_{\vM}$ as a function of $\vM = \vp n
+ \vm \sqrt{n}$.  In this notation, the desired bound we will
establish is
\begin{equation}
\label{symmetric}
\Delta_h \Delta_i  \Pr[I_{\vM} = j] = p_j \int_0^\infty\!\Delta_h \Delta_i \Delta_j F_t(\vM')\,dt = O(1 / n).
\end{equation}

Since we keep~$p_j$ fixed as $n \to \infty$, we see
from~(\ref{symmetric}) that we may treat~$h$, $i$, and~$j$
symmetrically.  Thus there are three cases to consider: (i)~$h$, $i$,
and~$j$ all distinct, (ii) $h = i \neq j$, (iii) $h = i = j$.  Let us
consider first the case~(i) of distinct $h, i, j$:
\begin{align*}
\Delta_h \Delta_i \Pr[I_{\vM} = j]
 &= \begin{aligned}[t]\int_0^\infty\!&\Pr[Y_{h, t} = M_h]  \Pr[Y_{i, t} = M_i] \\ &\times \Pr[Y_{j, t} = M_j - 1]\, p_j \left( \prod_{\ell \not\in \{h, i, j\}} \Pr[Y_{\ell, t} \geq M_{\ell}] \right)\,dt,\end{aligned}\\
\left| \Delta_h \Delta_i \Pr[I_{\vM} = j] \right|
 &\leq \int_0^\infty\!\frac{1}{\sqrt{2 \pi M_h}} \frac{1}{\sqrt{2 \pi M_i}} \Pr[Y_{j, t} = M_j - 1]\,p_j \,dt \\
 &= \frac{1}{\sqrt{2 \pi M_h}} \frac{1}{\sqrt{2 \pi M_i}}.
\end{align*}
If both $M_h$ and $M_i$ are
of order~$n$ (as when good strategies are used), then we obtain the
desired bound of $\left| \Delta_h \Delta_i \Pr[I_{\vM} = j] \right| = O(1 / n)$.
Otherwise, if (say) $M_i = o(n)$, we may instead use the bound
$$
\left| \Delta_i \Pr[I_{\vM} = j] \right| \leq \Pr[I_{\vm + e_i} \neq I_{\vm}] \leq \Pr[I_{\vM + e_i} = i],
$$
which is of much smaller order than $1 / n$.
Thus $\left| \Delta_h \Delta_i \Pr[I_{\vM} = j] \right| = O(1 / n)$
in any case.

Now suppose (ii) $h = i \neq j$.  In this case
\begin{align*}
\Delta_i \Delta_i \Pr[I_{\vM} = j]
 &= \begin{aligned}[t]\int_0^\infty\!&\big(\Pr[Y_{i, t} = M_i] - \Pr[Y_{i, t} = M_i +1]\big) \\ & \times \Pr[Y_{j, t} = M_j - 1]\,p_j \left( \prod_{\ell \not\in \{i, j\}} \Pr[Y_{\ell, t} \geq M_{\ell}] \right)\,dt,\end{aligned}\\
\left| \Delta_i \Delta_i \Pr[I_{\vM} = j] \right|
 &\leq \begin{aligned}[t]\int_0^\infty\!\big|&\Pr[Y_{i,t} = M_i] - \Pr[Y_{i,t} = M_i+1]\big| \\ & \times \Pr[Y_{j, t} = M_j - 1]\,p_j\,dt\end{aligned}\\
 &\leq \sup_t \big|\Pr[Y_{i, t} = M_i] - \Pr[Y_{i, t} = M_i + 1] \big|.
\end{align*}
Let $\lambda := p_i t$ and $M := M_i+1$.  For fixed $M$, the expression
\begin{align*}
\left| e^{- \lambda} \frac{\lambda^{M - 1}}{(M - 1)!} - e^{- \lambda}
        \frac{\lambda^M}{M!} \right|
  &=   e^{- \lambda} \frac{\lambda^{M - 1}}{(M - 1)!} \left| 1 -
        \frac{\lambda}{M} \right|
\end{align*}
is easily seen to be maximized when~$\lambda$ is one of the values $\lambda = M \mp \sqrt{M}$, and so
$$
\left| \Delta_i \Delta_i \Pr[I_{\vM} = j] \right| \leq \frac{1}{\sqrt{2 \pi M_i}} \frac{1}{\sqrt{M_i + 1}}.
$$
As in case~(i), we conclude that $\left| \Delta_i \Delta_i \Pr[I_{\vM} = j] \right| = O(1 / n)$.

In case~(iii) $h = i = j$, the desired bound follows from the equality
$$\sum_j \Delta_i\Delta_i\Pr[I_{\vM}=j]=\Delta_i\Delta_i 1 = 0$$
and the bounds for case (ii) $h = i \neq j$.
\end{proof}

Define the \emph{$\delta$-undercut} of~$\vx$ to be the mixed strategy
that increments~$\vx$ by $(k - 1) \delta$ in a uniformly random coordinate,
and decrements~$\vx$ by~$\delta$ in the remaining coordinates.
(In the discrete setting, $\delta$ will be a multiple of $1 / \sqrt{n}$.)
The preceding lemma implies
\begin{corollary} \label{cor:undercut}
  Uniformly in allocations~$\vx$ and~$\vy$ such that
  $|x_i - y_i| > (k - 1) \delta$ for each coordinate~$i$,
$$
\left| \Pr[\text{\rm $\delta$-undercut of~$\vx$ beats~$\vy$}] -
\Pr[\text{\rm $\vx$ beats~$\vy$}] \right| = O(\delta^2).
$$
\end{corollary}
\begin{proof}
  We prove this in the discrete setting; the result for the continuous
  setting follows by taking limits.  We assume without loss of
  generality that the bins are numbered so that $x_i>y_i$ for $1\leq
  i\leq \ell$ and $x_i<y_i$ for $\ell+1\leq i\leq k$.  Recall that
  $\vm=\max(\vx,\vy)$.  Let $\vx'$ denote the random $\delta$-undercut
  of $\vx$, and $\vm'=\max(\vx',\vy)$.  We have $$\Pr[\text{$\vx$
    beats $\vy$}] = \sum_{j=\ell+1}^k \Pr[I_{\vm}=j],$$
  and because
  $\vx$ and $\vy$ differ by more than $(k-1)\delta$ in each
  coordinate, $$\Pr[\text{$\vx'$ beats $\vy$}] = \sum_{j=\ell+1}^k
  \Pr[I_{\vm'}=j].$$
  For each coordinate $i$, $\E[m'_i]=m_i$, and
  $|m'_i-m_i|=O(\delta)$, so the corollary follows from the
  local-flatness property we proved in Lemma~\ref{L:flat}.
\end{proof}

\begin{lemma} \label{L:del}
  For any fixed~$\vp$ with $k\geq 3$ bins, for any coordinate~$i$, for
  any $\delta \geq 1 / \sqrt{n}$, and for any optimal
  strategy~$\alpha$, if~$\vx$ and~$\vy$ are independent draws
  from~$\alpha$, then $\Pr[|x_i - y_i| \leq \delta] = O(\delta)$,
  where the constant implicit in $O(\delta)$ depends only upon~$\vp$.
\end{lemma}

The same proof works for both the discrete and continuous versions of
the game.

\begin{proof}
  We construct a strategy $\beta$ which attempts to beat strategy
  $\alpha$ by undercutting it, to wit: $\beta$ picks an
  allocation~$\vx$ from~$\alpha$, but, rather than playing~$\vx$,
  strategy~$\beta$ instead plays the $\delta$-undercut of $\vx$.  By
  analyzing how~$\beta$ fares against~$\alpha$ we will be able to
  bound $\Pr[|x_i - y_i| \leq \delta]$.
  
  When~$\beta$ and~$\alpha$ are pitted against each other, we will
  take the viewpoint that~$\vx$ and~$\vy$ are independently drawn
  from~$\alpha$ and a fair coin is flipped; if the coin lands heads,
  then~$\alpha$ plays~$\vx$ and~$\beta$ plays the $\delta$-undercut
  of~$\vy$, while if the coin lands tails, then~$\alpha$ plays~$\vy$
  and~$\beta$ plays the $\delta$-undercut of~$\vx$.  Recall that $I$
  is the last bin to be cleared when the allocations are~$\vx$
  and~$\vy$.  Without the undercutting, bin~$I$ is owned by~$\beta$
  with probability $1/2$.  (We say that a player ``owns'' the last bin
  to be cleared if he was the player who placed more chips in that
  bin, or lost the coin toss in the event of tie.)  Let~$I'$ be the
  bin last to be cleared \emph{with} the undercutting.
  
  Letting~$E$ be the event that $|x_i - y_i|\leq (k-1)\delta$ for some
  coordinate~$i$ and $E^c$ its complement, we may express
\begin{align}
\Pr[\text{$\alpha$ beats~$\beta$}]
&=    \Pr[\{\text{$\alpha$ beats~$\beta$}\} \cap E]
        + \Pr[\text{$\alpha$ beats~$\beta$}\,|\,E^c] \Pr[E^c] \notag \\
&\leq \Pr[\{\text{$\alpha$ beats~$\beta$}\} \cap E] + \left(
        \mbox{$\frac{1}{2}$} + O(\delta^2) \right) \Pr[E^c]\notag
\intertext{by Corollary~\ref{cor:undercut}.  The optimality of~$\alpha$
implies}
0 &\leq \Pr[\{\text{$\alpha$ beats~$\beta$}\} \cap E] - \mbox{$\frac{1}{2}$}
          \Pr[E] + O(\delta^2).
\label{eqn:opt}
\end{align}
Now, conditioning on~$\vx$ and~$\vy$,
\begin{align*}
\Pr[\text{$\alpha$ beats~$\beta$}\,|\,\vx, \vy]  
&= \Pr[\text{$\beta$ owns~$I'\neq I$}\,|\,\vx,\vy]
 + \Pr[\text{$\beta$ owns~$I' =   I$}\,|\,\vx,\vy] \\
&\leq \Pr[I' \neq
I\,|\,\vx, \vy] + \Pr[\text{$\beta$ owns~$I$}\,|\,\vx,\vy].
\end{align*}
By Lemma~\ref{L:ichanges}, $\Pr[I' \neq I\,|\,\vx, \vy] = O(\delta)$.
To compute $\Pr[\text{$\beta$ owns~$I$}\,|\,\vx, \vy]$, we condition on~$I$.
For example, conditionally given the event $\delta < |x_I - y_I| < (k - 1)
\delta$, the player using~$\beta$ owns~$I$ with probability $1 / 2$ 
if the ``overcut bin''
[the bin with allocation incremented by $(k - 1) \delta$] chosen is 
\emph{not} bin~$I$
and with
probability~$1$ if it is.  Thus, if $\delta < |x_I - y_I| < (k - 1) \delta$,
then
$$
\Pr[\text{$\beta$ owns~$I$}\,|\,\vx, \vy, I]
  = \mbox{$\frac12$} \left( 1 - \mbox{$\frac{1}{k}$} \right) +
    \mbox{$\frac{1}{k}$} = \mbox{$\frac12$} + \mbox{$\frac{1}{2 k}$}.
$$
The other entries in the following formula are computed similarly:
$$ \Pr[\text{$\beta$ owns~$I$}\,|\,\vx, \vy, I] =
   \begin{cases}
   \frac12                 & \mbox{if\ }|x_I - y_I| > (k - 1) \delta \\
   \frac12 + \frac{1}{4 k} & \mbox{if\ }|x_I - y_I| = (k - 1) \delta \\
   \frac12 + \frac{1}{2 k} & \mbox{if\ }\delta < |x_I - y_I|
                                               < (k - 1) \delta \\
   \frac14 + \frac{3}{4 k} & \mbox{if\ }|x_I - y_I| = \delta \\
             \frac{1}{k}   & \mbox{if\ }|x_I - y_I| < \delta.
   \end{cases}
$$
Thus, conditional on~$\vx$ and~$\vy$ but not~$I$,
\begin{multline*}
\Pr[\text{$\alpha$ beats~$\beta$}\,|\,\vx, \vy] - \frac12 = O(\delta) \\
  + \frac{1}{2 k}\Pr[\delta \lhe |x_I - y_I| \lhe (k - 1)\delta\,|\,\vx, \vy]
  - \frac{k - 2}{2 k}\Pr[|x_I - y_I| \lhe \delta\,|\,\vx, \vy],
\end{multline*}
and so unconditionally
\begin{multline*}
\Pr[\{\text{$\alpha$ beats~$\beta$}\} \cap E] - \frac12 \Pr[E] = 
O(\delta) \Pr[E] \\
  + \frac{1}{2 k}\Pr[\delta \lhe |x_I - y_I| \lhe (k - 1) \delta]
  - \frac{k - 2}{2 k} \Pr[|x_I - y_I| \lhe \delta].
\end{multline*}
Substituting this into~\eqref{eqn:opt}, and then rearranging,
\begin{align*}
\frac{k - 2}{2 k}\Pr[|x_I - y_I| \lhe \delta]
  &\leq O(\delta) \Pr[E] + O(\delta^2) + \frac{1}{2 k} \Pr[\delta \lhe |x_I -
          y_I| \lhe (k - 1) \delta], \\ (k - 1) \Pr[|x_I - y_I| \lhe \delta]
 & \leq O(\delta) \Pr[E] + O(\delta^2) + \Pr[|x_I-y_I|\lhe (k-1)\delta].
\end{align*}

Recall from Theorem~\ref{T:first} that we have $O(\sqrt{n})$ bounds on
the overplay or underplay of the optimal strategy~$\alpha$.  It
follows that there is a positive constant~$q$ (depending on~$\vp$)
such that $\Pr[I = i\,|\,\vx, \vy] \geq q$ for any coordinate~$i$ and
plays~$\vx$ and~$\vy$ that~$\alpha$ might make.  Recalling also
that~$E$ is the event that $|x_i - y_i| \leq (k - 1) \delta$ for some
coordinate~$i$, we see that
$$
\Pr[E] \leq \frac{1}{q} \Pr[|x_I - y_I| \leq (k - 1) \delta].
$$
Letting $r := k-1$, we have then that, for some~$c$ and all $\delta \geq 1 /
\sqrt{n}$,
$$
\Pr[|x_I - y_I| \lhe \delta]
  \leq 
c \delta^2 + \frac{1 + c
\delta}{r} \Pr[|x_I - y_I| \lhe r \delta].
$$
When $r = 1$ (i.e.,\ $k = 2$) this inequality is uninformative, but otherwise
we may iterate it to show, for $j \geq 1$ and $\delta \geq 1 / \sqrt{n}$,
$$
\frac{\Pr[|x_I - y_I| \lhe \delta]}{(1 + c \delta) \times \cdots \times(1 + c 
         r^{j - 1} \delta)}
  \leq 
         c [\delta^2 + \cdots + r^{j - 1} \delta^2] + \frac{1}{r^j} \Pr[|x_I
         - y_I| \lhe r^j \delta].
$$
We take $j = \lceil \log(1 / \delta) / \log r \rceil$ so that
$(1 + c \delta) \times \cdots \times(1 + c r^{j - 1} \delta) = O(1)$ and
both terms on the right-hand side are $O(\delta)$.  Thus
$$
\Pr[|x_I - y_I| \leq \delta] = O(\delta).
$$
Recalling again that $\Pr[|x_I - y_I| \leq \delta] \asymp \Pr[\cup_i \{|x_i -
y_i| \leq \delta\}]$ yields the lemma.
\end{proof}

\begin{theorem} \label{T:ac}
  For competitive continuous Knock~'em Down with $k \geq 3$ bins, each
  marginal distribution of any optimal strategy is absolutely
  continuous with respect to Lebesgue measure, and has a square-integrable
  density function.
\end{theorem}
\begin{proof}
  If~$X$ and~$Y$ are i.i.d.\ draws from the marginal distribution,
  then, by Lemma~\ref{L:del}, $\Pr[|X - Y| \leq \delta]= O(\delta)$.
  From this it is straightforward (see, e.g., \cite[Theorem~2.12]{MR1333890})
  to check that the distribution is absolutely
  continuous with respect to Lebesgue measure, and that the density
  function is square-integrable.
\end{proof}

Note that the support of \emph{any\/} measure with nonzero absolutely
continuous part has positive Lebesgue measure, and that \emph{any\/}
subset of the line with positive Lebesgue measure has positive
1-dimensional Hausdorff measure.  Thus we have
\begin{corollary}
  For $k \geq 3$ bins, the support of any optimal strategy for
  continuous Knock~'em Down has Hausdorff dimension at least~$1$.
\end{corollary}
It is natural to guess that the true Hausdorff dimension is $k-1$, and
indeed that optimal strategies are absolutely continuous with
respect to $(k-1)$-dimensional Lebesgue measure.

\section{Existence of an optimal continuous strategy}
\label{S:optexists}

While every game in which each player has a finite number of options
will have a value (which of course is~$0$ for $n$-token
competitive Knock~'em Down), there are continuous games without a
value~\cite{MR20:265}.  Thus Theorem~\ref{T:optexists} below has
nontrivial content.

Recall the definition of~$K$ at~\eqref{Kdef}.  For the $n$-token game,
we regard a mixed strategy~$\alpha_n$ as a probability measure on
$k$-tuples $\vc = (x_1, \dots, x_k)$ with vanishing sum, as described
in Section~\ref{S:intro}, and we define the payoff function~$K_n$ on
this $\vc$-scale.  We say that a sequence $(\alpha_n)$ of strategies
for the $n$-token competitive Knock~'em Down games is
\emph{asymptotically optimal} if
$$
\min_{\vy\in A_n} K_n(\alpha_n, \vy) \to 0\mbox{\ \,as $n \to \infty$}.
$$
Here the $\min$ is taken over the finite (but growing, as $n \to
\infty$) set~$A_n$ of (normalized) actions (allocations) available 
for $n$-token Knock~'em Down.

\begin{theorem} \label{T:optexists}
  The continuous game has value~$0$, and there is at least one optimal
  strategy.  Indeed, any subsequential weak limit of any
  asymptotically optimal sequence $(\alpha_n)$ of strategies for
  $n$-token competitive Knock~'em Down is an optimal strategy for the
  continuous game.
\end{theorem}
\noindent
Key ingredients to the proof of this theorem are Theorem~\ref{T:first}
and Lemma~\ref{L:del}.  The converse to this theorem is proved in
Corollary~\ref{C:comp0}.

\begin{proof}[Proof of Theorem~\ref{T:optexists}]
  By Theorem~\ref{T:first} and Lemma~\ref{L:overplay}, any sequence
  $(\alpha_n)$ of asymptotically optimal strategies
  (probability measures) is tight.  Therefore
  there is a subsequential weak limit $\alpha$ which is a probability
  measure (see e.g.\ \cite[Theorem~3.1.9]{MR1267569}).  It is easy to
  check that~$\alpha$ is concentrated on $k$-tuples with vanishing
  sum, so that~$\alpha$ is a mixed strategy for the continuous game,
  and (using Lemma~\ref{L:del}) that~$\alpha$ has atomless marginals.
  We claim that~$\alpha$ is optimal, and then we see immediately that
  the continuous game has value~$0$.
  
  The proof that~$\alpha$ is optimal is rather routine.  To avoid
  double subscripts, we henceforth innocuously assume that the full
  sequence $(\alpha_n)$ converges weakly to~$\alpha$.  Fix any pure
  strategy~$\vy$ for the continuous game.  Let~$\vy_n$ denote a
  rounding of~$\vy$ to a pure strategy for~$K_n$; the details of the
  rounding procedure are irrelevant for our purposes here, as long as
  $\vy_n \to \vy$.  Using the facts that~$\alpha$ has atomless
  marginals, $\alpha_n \wto \alpha$, $K$ is continuous away from pairs
  $(\vx,\vy)$ for which $x_i=y_i$ for some $i$, and (for any $\eps>0$)
  $\left| K_n(\vx, \vy_n) - K(\vx, \vy) \right| \to 0$ uniformly for
  $\vx$ satisfying $|x_i - y_i| > \eps$ for $i = 1, \dots, k$, it
  follows (we omit the details) that $K_n(\alpha_n, \vy_n) \to
  K(\alpha, \vy)$.  But by the optimality of $\alpha_n$,
  $K_n(\alpha_n,\vy_n) \geq 0$ for every~$n$, so $K(\alpha,\vy)\geq 0$
  and~$\alpha$ is optimal.
\end{proof}

\begin{remark}
\label{R:unique}
(a)~If it happens to be true that there exists a \emph{unique}
optimal strategy~$\alpha_0$ for~$K$, then Theorem~\ref{T:optexists} 
implies that
\begin{equation}
\label{aconv}
\alpha_n \wto \alpha_0\mbox{\ \ as $n \to
\infty$}.
\end{equation}

(b)~If $p_1 = \cdots = p_k$, it is possible
to choose the strategies $\alpha_n$ to be symmetric.  If so, and if it happens
that there exists a unique symmetric optimal strategy~$\alpha_0$ for~$K$,
then~\eqref{aconv} holds.
\end{remark}

\section{Asymptotically optimal play of Knock~'em Down}
\label{S:asyopt}

The main result of this section (see Corollary~\ref{C:comp0}) is that
when an optimal strategy~$\alpha_0$ for the continuous game~$K$ is
``rounded'' to produce allocations for $n$-token Knock~'em Down, the
result is an asymptotically optimal strategy.  The drawback here is
that we do not know how to construct such an~$\alpha_0$, but we might
at least hope to find a not unreasonably suboptimal~$\apxopt$
for~$K$, such as the one discussed in Remark~\ref{approx-opt}.  We are
thus motivated to show more generally (see Theorem~\ref{T:comp1}) that
``rounding'' of any strategy~$\apxopt$ with atomless marginals and
bounded support gives a strategy for Knock~'em Down whose worst-case
payoff (i.e., payoff against an opponent playing the best possible
response) is asymptotically at least as large as the worst-case payoff
in game~$K$ from use of~$\apxopt$.

\begin{theorem} \label{T:comp1}
  If $\alpha_n \wto \apxopt$, where~$\apxopt$ has atomless marginals
  and bounded support, then
\begin{equation}
\label{comp1}
\liminf_{n \to \infty} \min_{\vy \in A_n} K_n(\alpha_n, \vy)
   \geq \inf_{\vy} K(\apxopt, \vy).
\end{equation}
\end{theorem}

\begin{proof}
Let us define
\begin{equation}
\label{worstpay}
\kappa_n := \min_{\vy \in A_n} K_n(\alpha_n, \vy),
\end{equation}
and let $\vy_n\in A_n$ be a strategy achieving this minimum.  Since
the sequence $(\alpha_n)$ converges weakly, it is tight, so
by Lemma~\ref{L:overplay}, $\vy_n$ remains bounded.
Let $n_{\ell} \uparrow \infty$ be any sequence for which
$\lim_{\ell\to\infty} \kappa_{n_\ell} = \liminf_{n\to\infty} \kappa_n$.
By compactness, we know that
there is a subsequence $\nt_{\ell} \uparrow \infty$ and a continuous-game
allocation~$\vy$ such that $\vy_{\nt_{\ell}} \to \vy$.

The supremum of
$|K_n(\vc,\vy)-K(\vc,\vy)|$ over any bounded set of
$(\vc,\vy)$ tends to $0$ as $n\to\infty$; in particular, ties do not
cause a problem here.  Since the sequence $(\alpha_n)$ is tight and
the sequence $\vy_n$ is bounded, we have that~$\kappa_n$ differs by
$o(1)$ from
$$
  \kappah_n := K(\alpha_n, \vy_n).
$$

For $\eps>0$ define
$$
F_\eps := \{ \vx:|x_i - y_i| \leq \eps\text{\ for some $i$}\},
$$
and let $F_\eps^c$ denote the complement of this set.  We have that
$$
K(\vc, \vy_{\nt_{\ell}}) \to K(\vc, \vy)\text{\ \ as $\ell \to \infty$}
$$
uniformly for~$\vc\in F_\eps^c$.  Thus, fixing $\eps > 0$,
\begin{equation*}
\kappah_{\nt_{\ell}} \geq
 \int_{F_\eps^c}\!K(\vc, \vy)\,\alpha_{\nt_{\ell}}(d \vc) -
\mbox{$\frac{1}{2}$} \alpha_{\nt_{\ell}}(F_\eps) - o(1).
\end{equation*}
Now, by an argument (omitted here) very much like one used in the detailed proof of
Theorem~\ref{T:optexists} one finds
$$
\lim_{\ell \to \infty} \kappah_{\nt_{\ell}}
 \geq K(\apxopt, \vy) - 2 \apxopt(F_\eps).
$$
Letting $\eps \downarrow 0$ and using the continuity of the marginals
of~$\apxopt$, we obtain
\begin{align*}
\liminf_{n\to\infty} \kappa_n = 
\liminf_{n\to\infty} \kappah_n = 
\lim_{\ell \to \infty} \kappah_{\nt_{\ell}} \geq
 K(\apxopt, \vy) &\geq
\inf_{\vy} K(\apxopt,\vy). \qedhere
\end{align*}
\end{proof}

\begin{remark}\label{approx-opt}
  Even very simple-minded mixed strategies can offer substantial
  improvement over pure strategies.  To illustrate this, we consider
  $k = 3$ and $\vp = (1 / 3, 1 / 3, 1 / 3)$.  Numerical
  experimentation suggested that the uniform distribution~$\apxopt$
  over the simplex of all triples $\vx = (x_1, x_2, x_3)$ summing to
  zero and satisfying $x_i \geq - 1 / 6$ for all~$i$ might be a good
  strategy for the continuous game.  Indeed, numerical explorations
  indicate that the best response to~$\apxopt$ is the pure strategy
  $(0, 0, 0)$ and thence that
$$
\mbox{RHS}\eqref{comp1} \doteq - 0.0101219
$$
[a huge improvement on the worst-payoff value $ - 1 / 6 \doteq
-0.1666667$ resulting from use of the na\"ive pure strategy $(0, 0,
0)$ in place of~$\apxopt$].  The strategy~$\apxopt$ was ``rounded'' to
a strategy~$\alpha_{180}$ for $180$-chip Knock~'em Down by taking
$\alpha_{180}$ to be uniform over allocations placing at least $58$
chips in each bin.  [Note $60 - (1/6) \sqrt{180} \doteq 57.8$.]  Then,
according to {\tt Mathematica}, the value of~$\kappa_{180}$
at~(\ref{worstpay}) is
$$
\kappa_{180} \doteq - 0.0165257,
$$
improving on the worst possible payoff of $ \doteq - 0.0920653$ for
the na\"ive allocation $(60, 60, 60)$.
\end{remark}

\begin{corollary}
\label{C:comp0} If $\alpha_n \wto \alpha_0$,
where~$\alpha_0$ is optimal for the continuous game~$K$,
then~$\alpha_n$ is asymptotically optimal for the $n$-token Knock~'em Down
game~$K_n$, in the sense that~$\kappa_n$ defined at~\eqref{worstpay} 
vanishes in
the limit.
\end{corollary}

\section{Open problems}
\label{S:open}

We have proved the existence of an optimal strategy for two-player
continuous Knock~'em Down, but we do not have an explicit description
of optimal play even when $k = 3$ and $p_1 = p_2 = p_3 = 1 / 3$.  We
know that the marginal distributions of optimal play are absolutely
continuous with respect to Lebesgue measure.  Consequently the set of
pure strategies supporting optimal play will have dimension at
least~$1$; perhaps the dimension is $k - 1$.

\section*{Acknowledgments}
We are grateful to Yuval Peres and Elchanan Mossel for useful
discussions, and the referee for helpful comments.  We thank Ed
Scheinerman, whose son Jonah played Knock~'em Down in his third-grade
math class taught by Bonnie Nagel, for bringing the game to our
attention.


\begin{thebibliography}{BFH01}

\bibitem[BF99]{BF_UMAP}
Arthur~T. Benjamin and Matthew~T. Fluet.
\newblock The best way to {Knock~'m Down}.
\newblock {\em The UMAP Journal}, 20(1):11--20, 1999.

\bibitem[BF00]{MR1767066}
Arthur~T. Benjamin and Matthew~T. Fluet.
\newblock What's best?
\newblock {\em Amer.\ Math.\ Monthly}, 107(6):560--562, 2000.
\MR{1767066}

\bibitem[BFH01]{MR2002g:91048}
Arthur~T. Benjamin, Matthew~T. Fluet, and Mark~L. Huber.
\newblock Optimal token allocations in {S}olitaire {Knock~'m Down}.
\newblock {\em Electron.\ J.\ Combin.}, 8(2):Research Paper~2, 8~pp., 2001.
\newblock In honor of Aviezri Fraenkel on the occasion of his 70th birthday.
\MR{1853253}

\bibitem[Fan53]{MR14:1109f}
Ky~Fan.
\newblock Minimax theorems.
\newblock {\em Proc.\ Nat.\ Acad.\ Sci.\ U.S.A.}, 39:42--47, 1953.
\MR{0055678}

\bibitem[Flu99]{Fluet}
Matthew~T. Fluet.
\newblock Searching for optimal strategies in {Knock~'m Down}, 1999.
\newblock Senior thesis, Harvey Mudd College, Claremont, CA.

\bibitem[FW06]{FW:solitaire}
James~A. Fill and David~B. Wilson.
\newblock Solitaire {Knock~'em Down}, 2006.
\newblock Manuscript.

\bibitem[Hun98]{Hunt}
Gordon Hunt.
\newblock Knock~'m down.
\newblock {\em Teaching Stat.}, 20(2):59--62, 1998.

\bibitem[Mat95]{MR1333890}
Pertti Mattila.
\newblock {\em Geometry of Sets and Measures in {E}uclidean Spaces, fractals
  and rectifiability}, volume~44 of {\em Cambridge Studies in Advanced
  Mathematics}.
\newblock Cambridge University Press, Cambridge, 1995.
\MR{1333890}

\bibitem[Str93]{MR1267569}
Daniel~W. Stroock.
\newblock {\em Probability Theory, an analytic view}.
\newblock Cambridge University Press, Cambridge, 1993.
\MR{1267569}

\bibitem[SW57]{MR20:265}
Maurice Sion and Philip Wolfe.
\newblock On a game without a value.
\newblock In {\em Contributions to the theory of games, vol.~3}, Annals of
  Mathematics Studies, no.~39, pages 299--306. Princeton University Press,
  Princeton, N. J., 1957.
\MR{0093742}

\end{thebibliography}

\end{document}